\documentclass[12pt,a4paper]{article}
\newtheorem{thm}{Theorem}
\newcommand{\la}{\lambda}
\newcommand{\cd}{\mathcal{D}}
\title{Franklin's argument proves an identity of Zagier}
\author{Robin Chapman\\
School of Mathematical Sciences\\ University of Exeter\\
Exeter, EX4 4QE, UK\\ \texttt{rjc@maths.ex.ac.uk}}
\date{4 September 2000}
\begin{document}
\maketitle

\begin{abstract}
Recently Zagier proved a remarkable $q$-series identity.
We show that this identity can also be proved by modifying
Franklin's classical proof of Euler's pentagonal number theorem.


\end{abstract}

\section{Introduction}

We use the standard $q$-series notation:
$$(a)_n=\prod_{k=1}^n(1-aq^{k-1})$$
where $n$ is a nonnegative integer or $n=\infty$.
Euler's pentagonal number theorem states that
\begin{equation}\label{pent}
(q)_\infty=1+\sum_{r=1}^\infty(-1)^r(q^{r(3r-1)/2}+q^{r(3r+1)/2}).
\end{equation}
Recently Zagier proved the following remarkable identity
\begin{thm}\label{Ztheorem}
\begin{equation}\label{Zag}
\sum_{n=0}^\infty[(q)_\infty-(q)_n]
=(q)_\infty\sum_{k=1}^\infty\frac{q^k}{1-q^k}
+\sum_{r=1}^\infty(-1)^r[(3r-1)q^{r(3r-1)/2}+3rq^{r(3r+1)/2}].
\end{equation}
\end{thm}
This is \cite[Theorem~2]{Z} slightly rephrased.

Equation (\ref{pent}) has a combinatorial interpretation.
The coefficient of $q^N$ in $(q)_\infty$ equals $d_e(N)-d_o(N)$
where $d_e(N)$ (respectively $d_o(N)$) is the number of partitions
of $N$ into an even (respectively odd) number of distinct parts.
Franklin \cite{F} showed that
$$d_e(N)-d_o(N)=\left\{
\begin{array}{cl}
(-1)^r&\textrm{if $N=\frac12r(3r\pm1)$ for a positive integer $r$,}\\
0&\textrm{otherwise.}
\end{array}\right.
$$
His proof was combinatorial. He set up what was almost an involution
on the set of partitions of $N$ into distinct parts. This ``involution''
reverses the parity of the number of parts. However there are certain
partitions for which his map is not defined. These exceptional partitions
occur precisely when $N=\frac12r(3r\pm1)$, and so account for the nonzero
terms on the right of~(\ref{pent}).

We show that Zagier's identity has a similar combinatorial interpretation,
which, miraculously, Franklin's argument proves at once. 

\section{Proof of Theorem~\ref{Ztheorem}}

We begin by recalling Franklin's ``involution''.
Let $\cd_N$ denote the set of partitions of $N$ into distinct parts and
let $\cd=\bigcup_{N=0}^\infty\cd_N$. For $\la\in\cd_N$ let $N_\la=N$,
$n_\la$ be the number of parts in $\la$ and $m_\la$ be the largest part of
$\la$ (if $\la$ is the empty partition of 0 let $m_\la=0$). Then
\begin{equation}\label{pentcomb}
(q)_\infty=\sum_{\la\in\cd}(-1)^{n_\la}q^{N_\la}.
\end{equation}

Let $\la$ be a non-empty partition in~$\cd$. Denote its smallest part by
$a_\la$. If the parts of $\la$ are $\la_1>\la_2>\la_3>\cdots$ let $b=b_\la$
denote the largest $b$ such that $\la_b=\la_1+1-b$ (so that $\la_k=\la_1+1-k$
if and only if $1\le k\le b$). If $\la\in\cd$ is not \emph{exceptional}
(we shall explain this term shortly), then we define a new partition $\la'$
as follows.
If $a_\la \le b_\la$ we obtain $\la'$ by removing the
smallest part from $\la$ and then adding 1 to the largest $a_\la$ parts of
this new partition.
If $a_\la>b_\la$ we obtain $\la'$ by subtracting 1 from the $b_\la$
largest parts of $\la$ and then appending a new part $b_\la$ to this
new partition.

The exceptional partitions are those for which this procedure breaks
down. We regard the empty partition as exceptional, also we regard those
for which $n_\la=b_\la$ and $a_\la=b_\la$ or $b_\la+1$. If $\la$ is not
exceptional, then neither is $\la'$ and $\la''=\la$
and $(-1)^{n_{\la'}}=-(-1)^{n_\la}$.
Thus on the right side of (\ref{pentcomb}) the contributions from
non-exceptional partitions cancel. The non-empty exceptional partitions
are of two forms: for each positive integer $r$ we have
$\la=(2r-1,2r-2,\ldots,r+1,r)$ for which $n_\la=r$, $m_\la=2r-1$ and
$N_\la=\frac12r(3r-1)$, and we have
$\la=(2r,2r-1,\ldots,r+2,r+1)$ for which $n_\la=r$, $m_\la=2r$ and
$N_\la=\frac12r(3r+1)$. Thus from (\ref{pentcomb}) we deduce
(\ref{pent}).

If $\la\in\cd$ is non-exceptional, then either $n_{\la'}=n_\la-1$
in which case $m_{\la'}=m_\la+1$ or $n_\la=n_\la+1$
in which case $m_{\la'}=m_\la-1$. In each case
$m_{\la'}+n_{\la'}=m_\la+n_\la$. It follows that in the sum
$$\sum_{\la\in\cd}(-1)^{n_\la}(m_\la+n_\la)q^{N_\la}$$
the terms corresponding to non-exceptional $\la$ cancel and
so we get only the contribution from exceptional~$\la$. Thus
\begin{equation}\label{Zagcomb}
\sum_{\la\in\cd}(-1)^{n_\la}(m_\la+n_\la)q^{N_\la}
=\sum_{r=1}^\infty(-1)^r[(3r-1)q^{r(3r-1)/2}+3rq^{r(3r+1)/2}].
\end{equation}
This sum occurs in (\ref{Zag}), which will follow
by analysing the left side of (\ref{Zagcomb}).

We break this into two sums. The first
$$\sum_{\la\in\cd}(-1)^{n_\la}m_\la q^{N_\la}$$
is dealt with in \cite[Theorem 5.2]{AJO}. We repeat their argument.
The coefficient of $q^N$ in $(q)_\infty-(q)_n$ is the sum of $(-1)^{n_\la}$
over all $\la\in\cd_N$ having a part strictly greater than~$n$.
Such a $\la$ is counted for exactly $m_\la$ different $n$ so that
\begin{equation}\label{msum}
\sum_{n=0}^\infty[(q)_\infty-(q)_n]
=\sum_{\la\in\cd}(-1)^{n_\la}m_\la q^{N_\la}.
\end{equation}

For each positive integer~$k$,
$$\frac{-q^k}{1-q^k}(q)_\infty
=(1-q)(1-q^2)\cdots(1-q^{k-1})(-q^k)(1-q^{k+1})\cdots.$$
The coefficient of $q^N$ in this product is the sum of $(-1)^{n_\la}$
over all $\la\in\cd_N$ having $k$ as a part. Such a $\la$ occurs for $n_\la$
distinct~$k$, and summing we conclude that
\begin{equation}\label{nsum}
-(q)_\infty\sum_{k=1}^\infty\frac{q^k}{1-q^k}
=\sum_{\la\in\cd}(-1)^{n_\la}n_\la q^{N_\la}.
\end{equation}
Combining (\ref{Zagcomb}), (\ref{msum}) and (\ref{nsum}) gives (\ref{Zag}).

\section{Another identity}

We can prove another identity by the same method. As before
Franklin's argument proves that
\begin{equation}\label{xeq}
\sum_{\la\in\cd}(-1)^{n_\la}x^{m_\la+n_\la}q^{N_\la}
=1+\sum_{r=1}^\infty(-1)^r[x^{3r-1}q^{r(3r-1)/2}+x^{3r}q^{r(3r+1)/2}].
\end{equation}

We now give an alternative way of expressing the left side of~(\ref{xeq}).
Consider the contribution due to the partitions
$\la$ with $m_\la$ equalling a fixed~$m\ge1$.
The product
$$(xq)_m=\prod_{j=1}^m(1-xq^j)
=\sum_{\la\in\cd\atop m_\la\le m}(-1)^{n_\la}x^{n_\la}q^{N_\la}$$
and so
$$\sum_{\la\in\cd\atop m_\la=m}(-1)^{n_\la}x^{m_\la+n_\la}q^{N_\la}
=x^m[(xq)_m-(xq)_{m-1}].$$
Hence
\begin{eqnarray*}
\sum_{\la\in\cd}(-1)^{n_\la}x^{m_\la+n_\la}q^{N_\la}
&=&1+\sum_{m=1}^\infty x^m[(xq)_m-(xq)_{m-1}]\\
&=&(1-x)\sum_{r=0}^\infty(xq)_mx^m\\
&=&\sum_{r=0}^\infty(x)_{r+1}x^r.
\end{eqnarray*}
Hence
\begin{equation}\label{S(x)}
\sum_{r=0}^\infty(x)_{r+1}x^r
=1+\sum_{r=1}^\infty(-1)^r[x^{3r-1}q^{r(3r-1)/2}+x^{3r}q^{r(3r+1)/2}].
\end{equation}
This can also be proved by noting that both sides satisfy the recurrence
$S(x)=1-qx^2-s^2x^3S(qx)$. This appears as exercise 10 in Chapter 2 of
\cite{A} whose solution is outlined in~\cite{Z}. Zagier deduces
(\ref{Zag}) from (\ref{S(x)}), essentially by carefully differentiating
with respect to $x$ and setting $x=1$.

\section{Acknowledgments}
The author wishes to thank George Andrews and Don Zagier for supplying
him with copies of \cite{AJO} and \cite{Z}.


\begin{thebibliography}{99}
\bibitem{A}
G.~E.~Andrews,
\emph{The Theory of Partitions},
Addison-Wesley, 1976 (reprinted Cambridge University Press 1998).
\bibitem{AJO}
G.~E.~Andrews, J.~Jim\'enez-Urroz, \& K.~Ono.
`Bizarre $q$-series identities and values of certain $L$-functions',
preprint.

\bibitem{F}
F.~Franklin, `Sur le d\'eveloppement du produit infini
$(1-x)(1-x^2)(1-x^3)(1-x^4)\ldots$',
\emph{C.\ R.\ Acad.\ Sci.\ Paris}, \textbf{92} (1881), 448--450.

\bibitem{Z}
D.~Zagier, `Vassiliev invariants and a strange identity related to the
Dedekind eta-function', \emph{Topology}, to appear.

\end{thebibliography}
\end{document}